\documentclass[11pt]{amsart} 
\usepackage{amscd,amssymb,amsxtra}
\usepackage{mathrsfs}  % \mathrsfs{A}
\DeclareSymbolFontAlphabet{\mathrsfs}{rsfs}
\usepackage[mathscr]{eucal} % \mathscr{A} 
\usepackage{comment}
\usepackage{xcolor}
\usepackage{enumitem} 
\usepackage[colorlinks,backref=page,citecolor=refkey]{hyperref}
\usepackage{citerange}
\usepackage{aliascnt}

\usepackage{marginnote}

\makeatletter
\let\@secnumfont\bfseries
\def\section{\@startsection{section}{1}%
  \z@{4\linespacing\@plus\linespacing}{\linespacing}%
  {\bfseries\centering}}
\def\introsection{\@startsection{section}{1}%
  \z@{3\linespacing\@plus\linespacing}{\linespacing}%
  {\bfseries\centering}}
\def\subsection{\@startsection{subsection}{2}%
   \z@{1.25\linespacing\@plus.7\linespacing}{.5\linespacing}%
   {\normalfont\bfseries}}
\def\subsectionsinline{\def\subsection{\@startsection{subsection}{2}%
  \z@{1\linespacing\@plus.7\linespacing}{-.5em}%
  {\normalfont\bfseries}}}

\makeatother

\numberwithin{equation}{section}

\newcommand{\mynewtheorem}[2]{
  \newaliascnt{#1}{equation}
  \newtheorem{#1}[#1]{#2}
  \aliascntresetthe{#1}
  \expandafter\def\csname #1autorefname\endcsname{#2}
}

\theoremstyle{definition}
\mynewtheorem{definition}{Definition}
\mynewtheorem{example}{Example}
\mynewtheorem{problem}{\color{blue}Problem}
\mynewtheorem{probsec}{\color{blue}Problem}
\mynewtheorem{exercise}{Exercise}
\mynewtheorem{question}{\color{blue}Question}
\mynewtheorem{project}{\color{blue}Project}
\mynewtheorem{construction}{Construction}
\mynewtheorem{task}{\color{blue}Task}
\mynewtheorem{otask}{\color{green}Obsolete Task}
\mynewtheorem{notation}{Notation}
\mynewtheorem{learn}{\color{blue}Learn}

\newtheorem*{definition*}{Definition}
\newtheorem*{example*}{Example}
\newtheorem*{problem*}{\color{blue}Problem}
\newtheorem*{probsec*}{\color{blue}Problem}
\newtheorem*{exercise*}{Exercise}
\newtheorem*{question*}{\color{blue}Question}
\newtheorem*{learn*}{\color{blue}Learn}
\newtheorem*{project*}{\color{blue}Project}
\newtheorem*{construction*}{Construction}
\newtheorem*{notation*}{Notation}

\theoremstyle{remark}
\mynewtheorem{note}{Note}
\mynewtheorem{remark}{Remark}
\mynewtheorem{data}{Data}

\newtheorem*{note*}{Note}
\newtheorem*{remark*}{Remark}
\newtheorem*{data*}{Data}

\theoremstyle{plain}
\mynewtheorem{theorem}{Theorem}
\mynewtheorem{corollary}{Corollary}
\mynewtheorem{lemma}{Lemma}
\mynewtheorem{proposition}{Proposition}
\mynewtheorem{conjecture}{Conjecture}
\mynewtheorem{claim}{Claim}
\mynewtheorem{proposal}{Proposal}
\mynewtheorem{conclusion}{Conclusion}
\mynewtheorem{hypothesis}{Hypothesis}
\mynewtheorem{assumption}{Assumption}

\newtheorem*{theorem*}{Theorem}
\newtheorem*{corollary*}{Corollary}
\newtheorem*{lemma*}{Lemma}
\newtheorem*{proposition*}{Proposition}
\newtheorem*{conjecture*}{Conjecture}
\newtheorem*{claim*}{Claim}
\newtheorem*{proposal*}{Proposal}
\newtheorem*{conclusion*}{Conclusion}
\newtheorem*{hypothesis*}{Hypothesis}
\newtheorem*{assumption*}{Assumption}

\newenvironment{proof*}[1][\proofname]{
  \begin{proof}[#1]}{  
\end{proof}}

\definecolor{refkey}{rgb}{0,.6,.4}

\renewcommand{\:}{\colon}
\newcommand{\CC}{{\mathbb C}}

\newcommand{\PP}{{\mathbb P}}

\newcommand{\RR}{{\mathbb R}}

\newcommand{\chiup}{\raise.5ex\hbox{$\chi$}}

\DeclareRobustCommand{\mstrut}{^{\vphantom{1*\prime y\vee M}}}

\newcommand{\temsquare}{\raise3.5pt\hbox{\boxed{ }}}

\DeclareFontFamily{U}{mathx}{}
\DeclareFontShape{U}{mathx}{m}{n}{<-> mathx10}{}
\DeclareSymbolFont{mathx}{U}{mathx}{m}{n}
\DeclareMathAccent{\widehat}{0}{mathx}{"70}
\DeclareMathAccent{\widecheck}{0}{mathx}{"71}
 
\DeclareMathSymbol{\bigtimes}{1}{mathx}{"91}

\newcommand{\upplus}{^{>0}}

\newcommand{\Rp}{\RR\upplus}

\definecolor{figblue}{rgb}{0.000,0.251,0.867}    % #0040DD
\definecolor{figred}{rgb}{0.843,0.000,0.082}     % #D70015
\definecolor{figgreen}{rgb}{0.082,0.502,0.239}   % #15803D
\definecolor{figpurple}{rgb}{0.608,0.176,0.557}  % #9B2D8E
\definecolor{figslate}{rgb}{0.216,0.255,0.318}   % #374151

\usepackage[all,2cell]{xy}
\usepackage{graphicx}\usepackage{epsf}

\usepackage{tikz}

\definecolor{refkey}{rgb}{0,.8,.2}\definecolor{labelkey}{rgb}{1,0,0} 

\DeclareMathOperator{\Bord}{Bord}
\DeclareMathOperator{\Imm}{Im}
\DeclareMathOperator{\Man}{Man}
\DeclareMathOperator{\Sym}{Sym}
\DeclareMathOperator{\Vect}{Vect}

\newcommand{\sF}{\mathscr{F}}
\newcommand{\sH}{\mathscr{H}}
\newcommand{\sU}{\mathscr{U}}
\newcommand{\sig}[2]{\leavevmode\hspace{3ex}\hbox to 40 pt{#1\hfil}{#2}\par}
\newcommand{\tVect}{t\!\Vect}
\newcommand{\tbord}{\Bord_{\langle n-1,n  \rangle}}
\renewcommand{\Im}{\Imm}

\usetikzlibrary{patterns,patterns.meta}
\usetikzlibrary{matrix,calc}
\usetikzlibrary{decorations.pathmorphing}

\tikzset{
  mfld/.style   = {black, line width=1.0pt},
  cover/.style  = {figblue, line width=0.8pt},
  nerve/.style  = {figgreen, line width=0.9pt},
  nvert/.style  = {figgreen, fill=figgreen},
  hatch/.style  = {pattern={Lines[angle=90,distance=7pt,line width=0.5pt]},
                   pattern color=figgreen},
  hatchfine/.style = {pattern={Lines[angle=90,distance=3.2pt,line width=0.4pt]},
                   pattern color=figgreen},
  uzero/.style  = {figblue, line width=0.9pt},
  uone/.style   = {figred,  line width=0.9pt},
  tube/.style   = {black, line width=0.7pt},
  xhatch/.style = {pattern={Lines[angle=55,distance=3.6pt,line width=0.4pt]},
                   pattern color=black!65},
}

\newcommand{\Mbig}{(5.50,2.80) ellipse [x radius=5.20, y radius=2.50]}
\newcommand{\setA}{(2.45,2.95) ellipse [x radius=1.85, y radius=1.20]}
\newcommand{\setB}{(4.05,3.65) ellipse [x radius=1.95, y radius=1.00]}
\newcommand{\setC}{(5.20,2.20) ellipse [x radius=3.45, y radius=0.95]}
\newcommand{\setD}{(7.45,3.35) ellipse [x radius=2.05, y radius=0.95, rotate=-22]}
\newcommand{\setE}{(8.15,2.05) ellipse [x radius=1.65, y radius=0.95]}

\newcommand{\Ucover}{%
  \draw[cover] \setA; \draw[cover] \setB; \draw[cover] \setC;
  \draw[cover] \setD; \draw[cover] \setE;
}

\newcommand{\tripleABC}{%
  \begin{scope}
    \clip \setA; \clip \setB; \clip \setC;
    \fill[hatchfine] (-1,-1) rectangle (12,7);
  \end{scope}}
\newcommand{\tripleCDE}{%
  \begin{scope}
    \clip \setC; \clip \setD; \clip \setE;
    \fill[hatchfine] (-1,-1) rectangle (12,7);
  \end{scope}}

\newcommand{\nervecoords}{%
  \coordinate (nA) at (0.00,1.10);
  \coordinate (nB) at (3.05,3.00);
  \coordinate (nC) at (3.95,0.05);
  \coordinate (nD) at (8.70,3.00);
  \coordinate (nE) at (8.95,0.55);
}
\newcommand{\covercoords}{%
  \coordinate (nA) at (2.45,2.95);
  \coordinate (nB) at (4.05,3.65);
  \coordinate (nC) at (5.20,2.20);
  \coordinate (nD) at (7.45,3.35);
  \coordinate (nE) at (8.15,2.05);
}
\newcommand{\nerveedges}{%
  \draw[nerve] (nA)--(nB) (nA)--(nC) (nB)--(nC) (nB)--(nD)
               (nC)--(nD) (nC)--(nE) (nD)--(nE);
}
\newcommand{\nerveverts}[1]{%
  \foreach \p in {nA,nB,nC,nD,nE} \fill[nvert] (\p) circle[radius=#1];
}

\newcommand{\Mturned}[2]{%   #1 = centre of the panel, #2 = scale
  \begin{scope}[shift={#1}]
    \begin{scope}[scale=#2]
      \begin{scope}[rotate=90]
        \begin{scope}[shift={(-5.50,-2.80)}]
          \draw[mfld] \Mbig;
          \Ucover
          \tripleCDE                      % U_0 n U_1 n U_2
          \node[black] at (7.45,3.80) {\small $U_0$};
          \node[black] at (9.05,1.75) {\small $U_1$};
          \node[black] at (4.60,2.20) {\small $U_2$};
        \end{scope}
      \end{scope}
    \end{scope}
  \end{scope}}
\newcommand{\cylend}[1]{%
  \begin{scope}[shift={(#1,0)}]
    \begin{scope}                       % the lens itself
      \clip (0,5.50) ellipse [x radius=0.44, y radius=1.30];
      \clip (0,4.10) ellipse [x radius=0.44, y radius=1.30];
      \fill[black!8] (-1,4) rectangle (1,6);
      \fill[xhatch]  (-1,4) rectangle (1,6);
    \end{scope}
    \begin{scope}                       % lower boundary: U_0's arc, blue
      \clip (-0.6,4.20) rectangle (0.6,4.80);
      \draw[uzero] (0,5.50) ellipse [x radius=0.44, y radius=1.30];
    \end{scope}
    \begin{scope}                       % upper boundary: U_1's arc, red
      \clip (-0.6,4.80) rectangle (0.6,5.40);
      \draw[uone] (0,4.10) ellipse [x radius=0.44, y radius=1.30];
    \end{scope}
  \end{scope}}
\newcommand{\LUzero}{(2.60,5.50) ellipse [x radius=0.44, y radius=1.30]}
\newcommand{\LUone}{(2.60,4.10) ellipse [x radius=0.44, y radius=1.30]}
\newcommand{\xzero}{6.60}
\newcommand{\xone}{13.60}

\newcommand{\Mfibre}{%
  \begin{scope}                       % hatch the lens: two clips intersect
    \clip \LUzero; \clip \LUone;
    \fill[xhatch] (1,3) rectangle (4,7);
  \end{scope}
  \draw[uzero] \LUzero;
  \draw[uone]  \LUone;}

\newcommand{\Mfibrelabels}{%
  \tikzset{lead/.style={black, line width=0.3pt}}
  \draw[lead] (1.78,6.15) -- (2.22,6.15);
  \draw[lead] (1.78,4.80) -- (2.23,4.80);
  \draw[lead] (1.78,3.45) -- (2.22,3.45);
  \node[black, left] at (1.74,6.15) {\small $U_0$};
  \node[black, left] at (1.74,4.80) {\small $U_0\cap U_1$};
  \node[black, left] at (1.74,3.45) {\small $U_1$};}

\newcommand{\blowup}{%
  \fill[black!8] (\xzero,4.20) rectangle (\xone,5.40);
  \draw[tube] (\xzero,5.40) -- (\xone,5.40);
  \draw[tube] (\xzero,4.20) -- (\xone,4.20);
  \draw[uzero] (\xzero,5.50) ellipse [x radius=0.44, y radius=1.30];
  \draw[uone]  (\xone,4.10)  ellipse [x radius=0.44, y radius=1.30];
  \cylend{\xzero}
  \cylend{\xone}
  \node[black] at (\xzero,7.15) {\small $\{0\}\times U_0$};
  \node[black] at (10.10,4.80)  {\small $\Delta^1\times(U_0\cap U_1)$};
  \node[black] at (\xone,2.45)  {\small $\{1\}\times U_1$};}

\newcommand{\dyy}{0.22}
\newcommand{\simpgap}[3]{%   #1 left node, #2 right node, k = #3
  \foreach \i in {0,...,\the\numexpr#3+1\relax}{%      faces: leftward, solid
    \pgfmathsetmacro{\yy}{(2*\i-(2*#3+2)/2)*\dyy}%
    \draw[->] ([yshift=\yy cm]#2.west) -- ([yshift=\yy cm]#1.east);}%
  \foreach \i in {0,...,#3}{%                    degeneracies: rightward, dotted
    \pgfmathsetmacro{\yy}{(2*\i+1-(2*#3+2)/2)*\dyy}%
    \draw[->,densely dotted] ([yshift=\yy cm]#1.east) -- ([yshift=\yy cm]#2.west);}}

\newcommand{\tfig}[1]{%
  \ifcsname tfigure#1\endcsname
    \csname tfigure#1\endcsname
  \else
    \GenericError{}{segal.tikz: no figure \string\tfig#1 defined}{}{}%
  \fi}
\expandafter\def\csname tfigure4\endcsname{%
\begin{tikzpicture}[scale=0.95]
  \draw[mfld] \Mbig;
  \Ucover
  \node[black]     at (10.05,4.75) {$M$};
  \node[figblue] at ( 1.15,0.55) {$\{U_\alpha\}$};
\end{tikzpicture}
  \vskip -.5pc
  
}

\expandafter\def\csname tfigure5\endcsname{%
\begin{tikzpicture}
  \begin{scope}
    \nervecoords
    \fill[hatch] (nA)--(nB)--(nC)--cycle;
    \fill[hatch] (nC)--(nD)--(nE)--cycle;
    \nerveedges
    \nerveverts{2.4pt}
  \end{scope}

  \begin{scope}[shift={(9.80,-0.15)}]
    \begin{scope}[scale=0.60]
      \draw[mfld] \Mbig;
      \Ucover
      \tripleABC
      \tripleCDE
      \covercoords
      \nerveedges
      \nerveverts{2.0pt}
    \end{scope}
  \end{scope}
\end{tikzpicture}
}

\expandafter\def\csname tfigure6\endcsname{%
\scalebox{0.70}{%
\begin{tikzpicture}
  \tikzset{sep/.style={black, line width=0.4pt}}

  \Mturned{(1.85,7.15)}{0.62}

  \draw[sep] (4.00,10.45) -- (4.00,3.75);         % between M and the formulas
  \draw[sep] (4.00,3.75)  -- (15.20,3.75);        % between formulas and nerve

  \node[align=center] at (9.70,7.15)
    {$\rho_i \colon M \longrightarrow [0,1], \qquad i = 0,\dots,4$\\[2.0ex]
     $\rho_0+\rho_1+\rho_2+\rho_3+\rho_4 = 1$};

  \begin{scope}[shift={(5.60,0.25)}]
    \begin{scope}[scale=1.02]
      \nervecoords
      \fill[hatch] (nA)--(nB)--(nC)--cycle;
      \fill[hatch] (nC)--(nD)--(nE)--cycle;
      \nerveedges
      \nerveverts{2.4pt}
      \node[black, right]       at (nD) {$0$};
      \node[black, below right] at (nE) {$1$};
      \node[black, below left]  at (nC) {$2$};
    \end{scope}
  \end{scope}
\end{tikzpicture}}
}

\expandafter\def\csname tfigure7\endcsname{%
\scalebox{0.96}{%
\begin{tikzpicture}
  \tikzset{sep/.style={black, line width=0.4pt}}

  \Mfibre

  \Mfibrelabels
  \node[black] at (2.60,2.30) {$M$};   % figure 6 omits this

  \draw[sep] (3.40,7.60) -- (3.40,1.55);
  \draw[sep] (3.40,1.55) -- (15.20,1.55);

  \blowup

  \draw[nerve] (\xzero,0.75) -- (\xone,0.75);
  \fill[nvert] (\xzero,0.75) circle[radius=2.4pt];
  \fill[nvert] (\xone,0.75) circle[radius=2.4pt];
  \node[black, below] at (\xzero,0.55) {$0$};
  \node[black, below] at (\xone,0.55) {$1$};
\end{tikzpicture}}
}

\expandafter\def\csname tfigure8\endcsname{%
\begin{tikzpicture}[line cap=round]
  \def\R{2.20}          % the circle
  \def\ra{2.32}         % U_0, just outside
  \def\rb{2.08}         % U_1, just inside
  \def\rc{2.52}         % U_2, further out (right-hand circle only)
  \tikzset{
    arc/.style={line width=2.2pt},
    circ/.style={black, line width=0.5pt},
    ptr/.style={black, line width=0.45pt, ->},   % plain tip: no arrows.meta
  }

  \begin{scope}
    \draw[circ] (0,0) circle[radius=\R];
    \draw[arc, figblue] (-25:\ra) arc[start angle=-25, end angle=205, radius=\ra];
    \draw[arc, figred]  (155:\rb) arc[start angle=155, end angle=370, radius=\rb];
    \node[figblue] at ( 55:2.95) {\small $U_0$};
    \node[figred]  at (295:2.95) {\small $U_1$};
    \node[black] at (0,0) {\small $U_0\cap U_1$};
    \draw[ptr] (-0.95,0.00) -- (-1.98,0.00);
    \draw[ptr] ( 0.95,-0.13) -- (1.96,-0.26);
  \end{scope}

  \begin{scope}[shift={(8.20,0)}]
    \draw[circ] (0,0) circle[radius=\R];
    \draw[arc, figblue]  ( 12:\ra) arc[start angle=12,  end angle=150, radius=\ra];
    \draw[arc, figgreen] (140:\rc) arc[start angle=140, end angle=235, radius=\rc];
    \draw[arc, figred]   (225:\rb) arc[start angle=225, end angle=382, radius=\rb];
    \node[figblue]  at ( 55:2.95) {\small $U_0$};
    \node[figred]   at (300:2.95) {\small $U_1$};
    \node[figgreen] at (180:3.20) {\small $U_2$};
  \end{scope}
\end{tikzpicture}
  \vskip -.5pc
  
}

\expandafter\def\csname tfigure9\endcsname{%
\scalebox{0.72}{%
\begin{tikzpicture}
  \tikzset{
    obj/.style = {circle, fill, inner sep=1.9pt},
    mor/.style = {black, line width=0.5pt, ->, shorten <=4pt, shorten >=4pt},
    big/.style = {black, line width=0.55pt, ->},
    squig/.style = {black, line width=0.5pt, ->, decorate,
      decoration={snake, amplitude=0.5mm, segment length=2.2mm,
                  post length=1.5mm}},
  }

  \draw[black, line width=0.8pt] (2.50,4.75) ellipse [x radius=2.40, y radius=1.44];
  \node[obj, black]   (cap) at (1.45,4.75) {};
  \node[obj, figblue] (uz)  at (3.15,5.40) {};
  \node[obj, figred]  (uo)  at (3.15,4.10) {};
  \draw[mor] (cap) -- (uz);
  \draw[mor] (cap) -- (uo);
  \node[black,   below] at (1.45,4.57) {\small $U_0\cap U_1$};
  \node[figblue, below] at (3.15,5.22) {\small $U_0$};
  \node[figred,  below] at (3.15,3.92) {\small $U_1$};
  \node[black] at (5.15,3.75) {${C}_{\sU}$};

  \draw[squig] (5.85,4.80) -- (7.85,4.80);

  \begin{scope}[shift={(7.95,0)}]
    \Mfibre
    \Mfibrelabels                               % no M label here
    \blowup
    \draw[big] (5.60,5.00) -- (3.60,5.00);      % blow-up covers M
    \draw[big] (10.10,3.40) -- (10.10,1.30);    % projection, centred on Delta^1

    \draw[nerve] (\xzero,0.75) -- (\xone,0.75);
    \draw[nerve, ->] (10.10,0.75) -- (8.30,0.75);
    \draw[nerve, ->] (10.10,0.75) -- (11.90,0.75);
    \foreach \x in {\xzero,10.10,\xone}
      \fill[nvert] (\x,0.75) circle[radius=2.4pt];
    \node[black, below] at (\xzero,0.55) {$0$};
    \node[black, below] at (10.10,0.55) {$01$};
    \node[black, below] at (\xone,0.55) {$1$};
  \end{scope}
\end{tikzpicture}}
  \vskip -.5pc
  
}

\expandafter\def\csname tfigure10\endcsname{%
\begin{tikzpicture}[scale=0.8,>=latex]
  \fill[figgreen!12] (-4.9,0) rectangle (4.9,-3.25);
  \draw[black, line width=0.8pt, dash pattern=on 6pt off 4pt, <->] (-4.9,0) -- (4.9,0);
  \draw[figblue, line width=1.0pt, ->] (0,-0.09) -- (0,-3.25);
  \draw[figblue, line width=0.8pt, fill=white] (0,0) circle[radius=0.09];
  \node[black]     at (4.45,-0.32) {$t$};
  \node[figblue]   at (0.32,-2.80) {$\tau$};
  \node[figgreen]  at (4.50,-2.80) {$z$};
\end{tikzpicture}}

\expandafter\def\csname tfigure11\endcsname{%
\begin{tikzpicture}[scale=0.7,>=latex]
  \tikzset{
    orient/.style = {figred,  line width=1.0pt, ->},
    tarrow/.style = {figblue, line width=1.0pt, ->},
    pt/.style     = {figgreen, fill=figgreen},   % objects (boundaries) are green
    ival/.style   = {black, line width=1.1pt},
    maps/.style   = {black, line width=0.6pt, |-to},
  }
  \begin{scope}[shift={(-1.2,0)}]   % left column, pulled away from the |-> arrows
  \draw[orient] (3.05,1.55) -- (3.907,1.55);
  \shade[left color=white, right color=black!50] (2.85,0.955) rectangle (3.45,1.045);
  \shade[left color=black!50, right color=white] (3.45,0.955) rectangle (4.05,1.045);
  \fill[pt] (3.45,1.0) circle[radius=3.2pt];
  \node[black] at (4.45,1.0) {$+$};
  \draw[tarrow] (2.75,0.45) -- (4.30,0.45);
  \shade[left color=white, right color=black!50] (0.75,-1.545) rectangle (1.35,-1.455);
  \shade[left color=black!50, right color=white] (6.15,-1.545) rectangle (6.75,-1.455);
  \draw[ival] (1.35,-1.5) -- (6.15,-1.5);
  \foreach \x in {1.35,6.15} \fill[pt] (\x,-1.5) circle[radius=3.2pt];
  \node[black] at (0.35,-1.5) {$+$};
  \node[black] at (7.15,-1.5) {$+$};
  \draw[orient] (0.85,-0.95) -- (1.707,-0.95);
  \draw[orient] (5.65,-0.95) -- (6.507,-0.95);
  \node[black] at (3.75,-0.95) {$\tau$};
  \draw[tarrow] (3.00,-2.05) -- (4.50,-2.05);
  \end{scope}
  \draw[maps] (8.30,1.0)  -- (10.30,1.0)  node[midway,above] {$F$};
  \draw[maps] (8.30,-1.5) -- (10.30,-1.5) node[midway,above] {$F$};
  \node[figgreen] at (14.00,1.0)  {$\sH$};
  \node[black] at (14.00,-1.5)
     {$\bigl({\color{figgreen}\sH}\xrightarrow{\;e^{-\tau H/\hbar}\;}{\color{figgreen}\sH}\bigr)$};
\end{tikzpicture}}

\newcommand{\germband}[8]{%
  \foreach \k in {15,14,...,0}{%
    \pgfmathsetmacro{\gt}{0.02*\k}%
    \pgfmathsetmacro{\gg}{45*(1-\k/15)}%
    \draw[black!\gg, line width=1.4pt]
      ({#1-\gt*#6},{#2-\gt*#7}) ellipse [x radius=#3, y radius=#4, rotate=#5];
    \ifnum#8=1
      \draw[black!\gg, line width=1.4pt]
        ({#1+\gt*#6},{#2+\gt*#7}) ellipse [x radius=#3, y radius=#4, rotate=#5];
    \fi}}

\expandafter\def\csname tfigure13\endcsname{%
\begin{tikzpicture}[scale=0.6,>=latex]
  \tikzset{
    bdy/.style    = {figgreen, line width=1.0pt},
    bord/.style   = {black,    line width=1.0pt},
    orient/.style = {figred,   line width=1.0pt, ->},
    maps/.style   = {black,    line width=0.6pt, |-to},
    hole/.style   = {black,    line width=0.8pt},
  }
  \begin{scope}[shift={(2.8,3.4)}]
    \germband{0}{0}{0.40}{0.95}{-18}{0.951}{-0.309}{1}
    \draw[bdy] (0,0) ellipse [x radius=0.40, y radius=0.95, rotate=-18];
    \draw[orient] (0.063,-0.021) -- (1.014,-0.330);   % perpendicular to Y
    \node[figgreen] at (-0.95,-1.40) {$Y^{n-1}$};
  \end{scope}
  \draw[bord] (0,0.85) .. controls (1.8,0.35) and (3.3,0.45) .. (4.97,1.32);
  \draw[bord] (5.53,-0.22) .. controls (4.6,-0.55) and (3.6,-1.1) .. (3.595,-2.15);
  \draw[bord] (0,-0.85) .. controls (1.2,-0.85) and (1.7,-1.0) .. (1.705,-1.95);
  \foreach \cx/\cy/\r in {2.05/0.12/10, 3.30/-0.62/-15}{
    \begin{scope}[shift={(\cx,\cy)}, rotate=\r, scale=1.35]
      \draw[hole] (-0.32,0.04) to[bend right=45] (0.32,0.04);
      \draw[hole] (-0.22,-0.02) to[bend left=45] (0.22,-0.02);
    \end{scope}}
  \node[black] at (2.00,-0.95) {$X$};
  \germband{0}{0}{0.35}{0.85}{0}{1}{0}{0}
  \germband{5.25}{0.55}{0.38}{0.82}{20}{-0.940}{-0.342}{0}
  \germband{2.65}{-2.05}{0.95}{0.26}{-6}{0}{1}{0}
  \draw[bdy] (0,0) ellipse [x radius=0.35, y radius=0.85];
  \draw[bdy] (5.25,0.55) ellipse [x radius=0.38, y radius=0.82, rotate=20];
  \draw[bdy] (2.65,-2.05) ellipse [x radius=0.95, y radius=0.26, rotate=-6];
  \node[figgreen] at (-1.00,0.15)  {$Y_1$};
  \node[figgreen] at (6.25,1.15)   {$Y_2$};
  \node[figgreen] at (2.65,-2.90)  {$Y_3$};
  \draw[orient] (0.017,0) -- (1.017,0);
  \draw[orient] (5.206,0.534) -- (4.266,0.192);
  \draw[orient] (2.65,-2.123) -- (2.65,-1.123);
  \draw[maps] (7.30,3.4) -- (9.30,3.4) node[midway,above] {$F$};
  \draw[maps] (7.30,-0.4) -- (9.30,-0.4) node[midway,above] {$F$};
  \node[figgreen, anchor=west, inner sep=0pt] at (12.486,3.4) {$F(Y)$};
  \node[black] at (15.6,-0.4) {$\Bigl(F(X)\colon
     {\color{figgreen}F(Y_1)}\otimes{\color{figgreen}F(Y_2)}\otimes{\color{figgreen}F(Y_3)}
     \longrightarrow{\color{figgreen}\mathbb{C}}\Bigr)$};
\end{tikzpicture}}

\newcommand{\tcechderham}{%
\begin{tikzpicture}[baseline=(current bounding box.center),>=latex]

\matrix (m) [
  matrix of math nodes,
  ampersand replacement=\&,   % REQUIRED: a \matrix inside a macro has its
  nodes={inner sep=2pt},
  execute at begin cell=\displaystyle,   % full-size \prod
  column sep=1.0cm,
  row sep=1.0cm,
]{
  \vdots \& \vdots \& \vdots \& \\
  \prod\Omega^{2}(U_{\alpha_0})
    \& \prod\Omega^{2}(U_{\alpha_0}\cap U_{\alpha_1})
    \& \prod\Omega^{2}(U_{\alpha_0}\cap U_{\alpha_1}\cap U_{\alpha_2}) \& \cdots \\
  \prod\Omega^{1}(U_{\alpha_0})
    \& \prod\Omega^{1}(U_{\alpha_0}\cap U_{\alpha_1})
    \& \prod\Omega^{1}(U_{\alpha_0}\cap U_{\alpha_1}\cap U_{\alpha_2}) \& \cdots \\
  \prod\Omega^{0}(U_{\alpha_0})
    \& \prod\Omega^{0}(U_{\alpha_0}\cap U_{\alpha_1})
    \& \prod\Omega^{0}(U_{\alpha_0}\cap U_{\alpha_1}\cap U_{\alpha_2}) \& \cdots \\
};

\coordinate (o)   at ($(m-4-1.south west)+(-0.55,-0.50)$);
\coordinate (top) at ($(o |- m-1-1.north)+(0,0.25)$);
\coordinate (rgt) at ($(m-2-4.east |- o)+(0.40,0)$);
\draw (o) -- (top);      % de Rham degree
\draw (o) -- (rgt);      % Cech degree

\node[left] at ($(m-2-1.west)+(-0.70,0)$) {$2$};
\node[left] at ($(m-3-1.west)+(-0.70,0)$) {$1$};
\node[left] at ($(m-4-1.west)+(-0.70,0)$) {$0$};

\node[below] at ($(m-4-1.south)+(0,-0.62)$) {$0$};
\node[below] at ($(m-4-2.south)+(0,-0.62)$) {$1$};
\node[below] at ($(m-4-3.south)+(0,-0.62)$) {$2$};

\begin{scope}[shorten <=1pt, shorten >=1pt]
  \draw[->] (m-4-1) -- (m-3-1) node[midway,left]  {$d$};
  \draw[->] (m-4-1) -- (m-4-2) node[midway,below] {$\delta$};
\end{scope}

\end{tikzpicture}
}

\newcommand{\tsimpobj}{%
\begin{tikzpicture}[baseline=(n0.base)]
  \node (n0) at (0,0)   {$NC^0$};
  \node (n1) at (2.6,0) {$NC^1$};
  \node (n2) at (5.2,0) {$NC^2$};
  \node (n3) at (7.6,0) {$\cdots$};
  \simpgap{n0}{n1}{0}   % 2 faces, 1 degeneracy
  \simpgap{n1}{n2}{1}   % 3 faces, 2 degeneracies
  \simpgap{n2}{n3}{2}   % 4 faces, 3 degeneracies
\end{tikzpicture}
}

\begin{document}

\abovedisplayskip18pt plus4.5pt minus9pt
\belowdisplayskip \abovedisplayskip
\abovedisplayshortskip0pt plus4.5pt
\belowdisplayshortskip10.5pt plus4.5pt minus6pt
\baselineskip=15 truept
\marginparwidth=55pt

\makeatletter
\renewcommand{\tocsection}[3]{%
  \indentlabel{\@ifempty{#2}{\hskip1.5em}{\ignorespaces#1 #2.\;\;}}#3}
\renewcommand{\tocsubsection}[3]{%
  \indentlabel{\@ifempty{#2}{\hskip 2.5em}{\hskip 2.5em\ignorespaces#1%
    #2.\;\;}}#3} 
\renewcommand{\tocsubsubsection}[3]{%
  \indentlabel{\@ifempty{#2}{\hskip 5.5em}{\hskip 5.5em\ignorespaces#1%
    #2.\;\;}}#3} 
%  \indentlabel{\hskip 4em#3}}
\def\@makefnmark{%
  \leavevmode
  \raise.9ex\hbox{\fontsize\sf@size\z@\normalfont\tiny\@thefnmark}} 
\def\multfoot{\textsuperscript{\tiny\color{red},}}
\def\footref#1{$\textsuperscript{\tiny\ref{#1}}$}
%             \footnote{Note\label{SYM}} ... in footnote~\footref{SYM}
\makeatother

\setcounter{tocdepth}{2}

%**end of header

% lasteq@  7
% lastsec@  4
% lastthm@  3
% lastfig@ 13

 \title[October 12, 2026]{On the work of Graeme Segal} %% replace \today with short title in final version

 \author[D. S. Freed]{Daniel S.~Freed}
 \address{Harvard University \\ Department of Mathematics \\ Science Center
Room 325 \\ 1 Oxford Street \\ Cambridge, MA 02138}
 \email{dafr@math.harvard.edu}

% \dedicatory{}
\date{October 12, 2026}
\maketitle

Graeme Segal was awarded the Chern Medal at the 2026 International Congress of
Mathematicians (ICM) ``for his visionary mathematical insights, that have had
an enduring influence in a wide range of fields, including topology,
mathematical physics, representation theory and category theory.''  In her
Citation for the Chern Medal~\cite{IMU}, Nathalie Wahl concludes:
 
        \smallskip
	\begin{quotation}\footnotesize
 Beyond technical achievements, Segal's impact lies in his compelling vision
and clarifying ideas. Many of his ideas have become so embedded into the
subject that their origins are often forgotten, yet several areas of
mathematics and mathematical physics are now unthinkable without them. His work
has integrated topology, geometry, representation theory, and physics into a
unified geometric perspective. He has illuminated difficult concepts with
remarkable simplicity, making them accessible to mathematicians and physicists
alike, and his generosity in sharing ideas has amplified his influence far
beyond his published papers. His work has not only solved deep problems but has
changed the way entire fields are conceptualized.
	\end{quotation}
        \smallskip

In this note, an expanded version of the Laudatio I gave at the ICM, I will
illustrate these features of Segal's mathematical work.  After a high-level
overview~(\S\ref{sec:1}), I turn to two specific episodes: classifying spaces
of categories~(\S\ref{sec:2}) and Segal's mathematical conceptualization of
quantum field theory~(\S\ref{sec:4}).

Segal began his graduate studies at Cambridge as a student of Sir William
Hodge.  After a year Segal moved to Oxford to work with (the much later Sir)
Michael Atiyah.  Still, in that initial year Hodge had great influence, perhaps
none as important as his first exhortation to a young Segal: Read Chern's
papers!  Thus began Segal's lifelong engagement with Chern's work.  While they
have different styles and approaches, both Chern and Segal are profound
geometers whose ideas have become rooted in the fabric of mathematics.  For
this and many more reasons, Segal is an ideal recipient of this year's Chern
Medal.

Warm thanks to Graeme for many recent discussions related to this Laudatio, for
providing his written and verbal perspectives (some quoted at length
in~\S\ref{sec:4}), and more generally for many decades of mathematics and
friendship.  Claude was used to help with manuscript preparation: figures,
references, and self-refereeing.

  \begin{figure}[htbp]
  \centering
  \includegraphics[scale=.5]{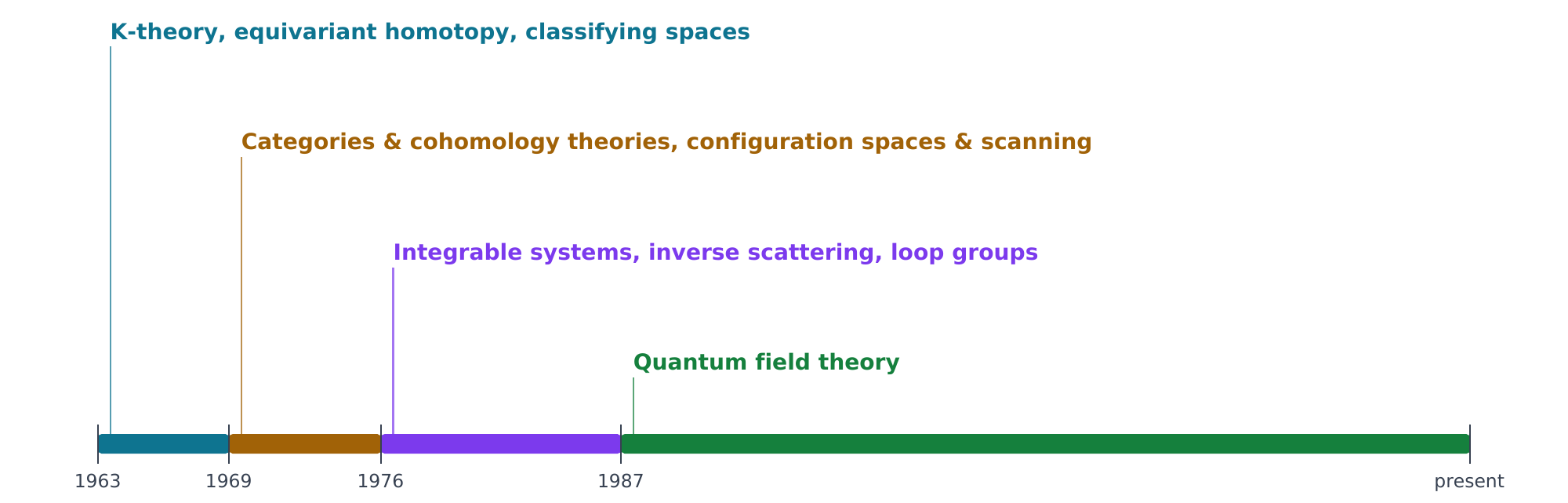}
  \vskip -.5pc
  \caption{Summary timeline of Segal's mathematical work}\label{fig:1}
  \end{figure}

   \section{Overview}\label{sec:1}
% lastsubsec@000

Segal has always maintained a broad vision and wide-ranging interest within
mathematics; his research at any given time reflects that breadth of approach.
Nevertheless, there is a clear arc to his work: it begins in algebraic
topology, proceeds to applications in geometry, continues to loop groups and
integrable systems, and settles into a focus on a mathematical
conceptualization of Wick-rotated quantum field theory.  See \autoref{fig:1}
for a timeline.

Along the way Segal proved many landmark theorems.  Here is a list of some of
his major results:

\medskip
\sig{1968}{Localization in equivariant $K$-theory and the Lefschetz fixed point theorem~\cite{AS1}}
\sig{1969}{Atiyah--Segal completion theorem in equivariant $K$-theory~\cite{AS2}}
\sig{1973}{Configuration spaces model iterated loop spaces~\cite{S2}}
\sig{1974}{Recognition of infinite loop spaces via $\Gamma $-spaces~\cite{S3}}
\sig{1976}{McDuff--Segal group-completion theorem~\cite{McS}}
\sig{1979}{Topology of spaces of rational functions; stability~\cite{S4}}
\sig{1981}{Positive-energy representations of loop groups and Borel--Weil~\cite{S5,PS}}
\sig{1985}{Segal--Wilson Grassmannian construction of KdV/KP solutions~\cite{SW}}
\sig{2001}{Open-closed 2d TQFT and D-branes as boundary conditions~\cite{MoS}}
\sig{2021}{Positivity of energy and Wick rotation~\cite{KS}}
\medskip 
 
\noindent 
 The early localization theorem in equivariant $K$-theory is part of Segal's
thesis~\cite{S12,S1} and appears jointly with Atiyah in the series of Annals papers
on the Atiyah--Singer index theorem.  The 2001~work on 2-dimensional
topological field theories is joint with Greg Moore.  The more recent paper on
positivity of energy in quantum field theory and Wick rotation, which proceeds
via a class of allowable complex metrics on a smooth manifold, is joint work
with Maxim Kontsevich.

One theme of this note is that Segal's influence goes far beyond these and
other particular theorems that he proved.  He also plants seeds that blossom in
subsequent work, opening up new research directions and even whole approaches
to existing fields.  Here is a sample:

\medskip
\sig{1968}{Classifying space of a topological category~\cite{S6}}
\sig{1968}{Equivariant $K$-theory as geometric representation theory~\cite{S1}}
\sig{1970}{Segal's Burnside ring conjecture~\cite{S7}}
\sig{1974}{$\Gamma $-spaces and the Segal condition~\cite{S3}}
\sig{1978}{Classifying spaces of foliations~\cite{S8}}
\sig{1979}{Segal's scanning technique~\cite{S4}}
\sig{1985}{Segal--Wilson Grassmannian~\cite{SW}}
\sig{1988}{Segal axioms for 2-dimensional conformal field theory~\cite{S9,S10}}
\sig{1988}{Geometric models for elliptic cohomology~\cite{S11}}
\sig{1995}{Floer homotopy theory~\cite{CJS}}
\medskip 
 
\noindent 
 The last item, Floer homotopy theory, was initiated jointly with Ralph Cohen
and John Jones.  It illustrates the observation that Segal's work is often a
point of departure for a subsequent flowering of new ideas: see~\cite{M1, KM,
M2, C1, C2, LS, M3, Kh, St, KLS1, KLS2, AB1, AB2}.  I will give a more detailed
account of other points of departure in the sequel.
 
Before going there, I share some thoughts on a question I was asked when
preparing this Laudatio: \emph{What kind of mathematician is Graeme Segal?}
Here, in bullet points, are my responses:

  \begin{itemize}
  \item Segal's motivation is to \emph{understand}; he is not a ``big game
hunter''\footnote{Segal's terminology}
  \item Segal is a master at conceptualizing structure
  \item Segal penetrates to the core of a mathematical idea
  \item Segal is often ahead of his time
  \item Segal inspires through his insightful exposition and communication
  \item Segal freely shares ideas and often leaves others to develop them
  \item Segal has incisive geometric intuition 
  \item Segal's mathematics is elegant and innovative
  \item Segal is a modest mathematician
  \item Segal is \emph{not} the type of mathematician who typically wins prizes
  \item Segal \emph{is} the type of mathematician who will be \emph{last} to be replaced by AI
  \end{itemize}

  \begin{figure}[htbp]
  \centering
  \includegraphics{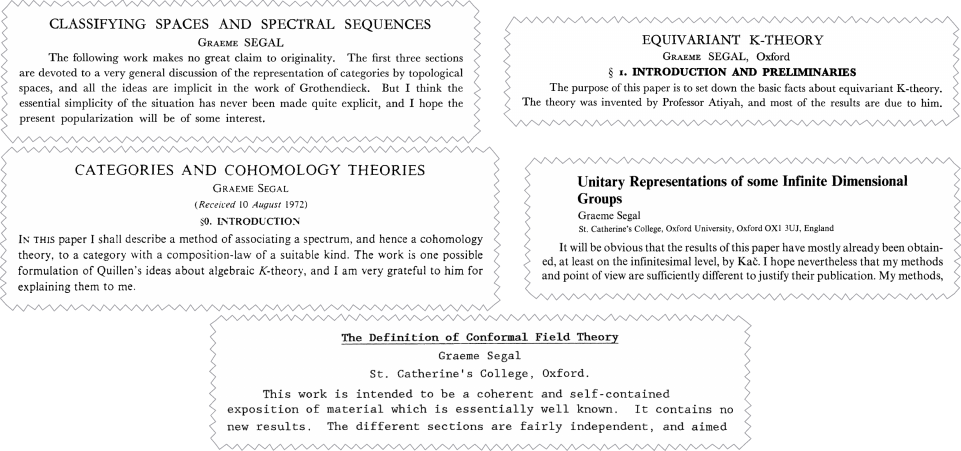}
  \vskip -.5pc
  \caption{Openings of five of Segal's most influential papers~\cite{S6,S1,S3,S5,S10}}\label{fig:2}
  \end{figure}

\noindent
 Segal's modesty is illustrated by the introductions to several of his papers;
see \autoref{fig:2}.  He either attributes the results to others, or he
disavows any originality in the results.  Yet these are among his most original
and influential works.

  \begin{figure}[htbp]
  \centering
  \includegraphics[scale=.9]{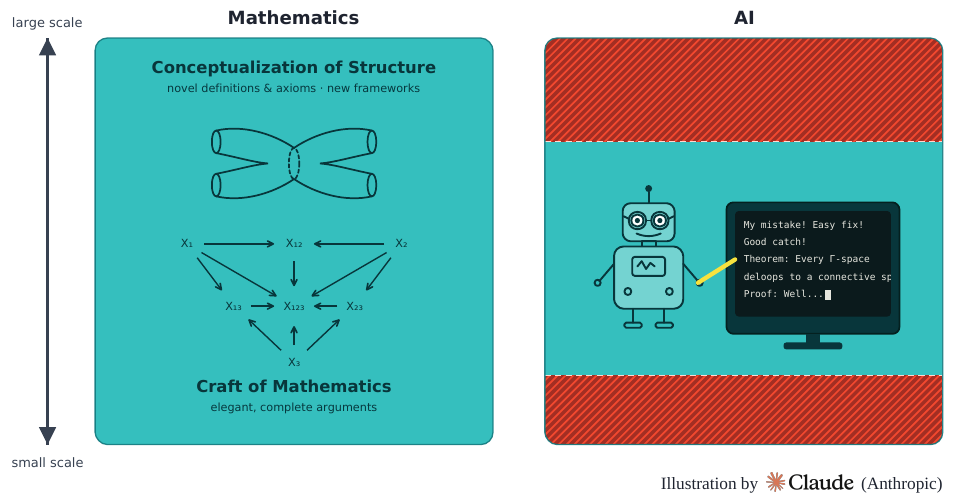}
  \vskip -.5pc
  \caption{Great mathematics vs.~AI}\label{fig:3}
  \end{figure}

By the last bullet point, which is the most topical answer to the question
\emph{What kind of mathematician is Graeme Segal?}, I mean something quite
specific.  Great mathematics and great mathematicians work at all scales, often
simultaneously.  At the smallest scale is the craft of mathematics: correct and
elegant reasoning compiled into insightful proofs that penetrate to the core.
At the largest scale are novel ideas that open up new research directions:
conceptualizations of structure common to many examples, and new definitions
and axioms that apply in unexpected ways.  By contrast, impressive as AI is at
present, it exhibits a cutoff at both extremes; see the schematic in
\autoref{fig:3}.  The situation with AI is changing rapidly, and it is unclear
where it will end up.  But that is not my subject here.  Rather, I will use two
episodes in Segal's mathematical output to demonstrate how he thrives at both
extremes.

   \section{Classifying spaces of categories}\label{sec:2}
% lastsubsec@000

\begin{figure}[htbp]
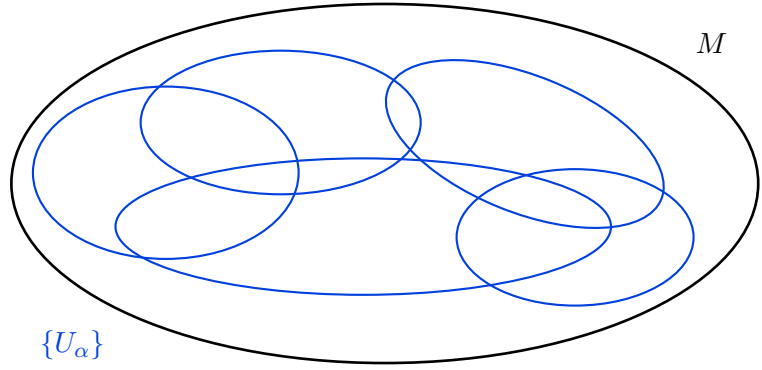

\centering
  \tfig4
  \vskip -.5pc
  \caption{An open cover $\{U_\alpha\}$ of $M$}\label{fig:4}
\end{figure}

As a graduate student Segal was musing on Andr\'e Weil's proof of the de Rham
theorem~\cite{W}.  Let $M$~be a smooth manifold, and suppose $\{U_\alpha
\}_{\alpha \in A}$ is an open covering of~$M$, as in~\autoref{fig:4}.  Weil's
argument proceeds via a double complex in which the horizontal differential
computes \v Cech cohomology and the vertical differential computes de Rham
cohomology:
  \begin{equation*}\label{eq:1}
\tcechderham
  \end{equation*}
(This double complex was later popularized by Bott--Tu~\cite{BT}.  It is worth
noting that Weil's paper has other useful double complexes, for example
versions that compute homology rather than cohomology, and he also addresses
Poincar\'e duality.)  As an early manifestation of his geometric bent, Segal
pondered the following question:
  $$ \textnormal{\emph{How can I realize Weil's \v Cech--de Rham double
     complex geometrically?}} $$%
That led to the following sequence of moves. 

\begin{figure}[htbp]
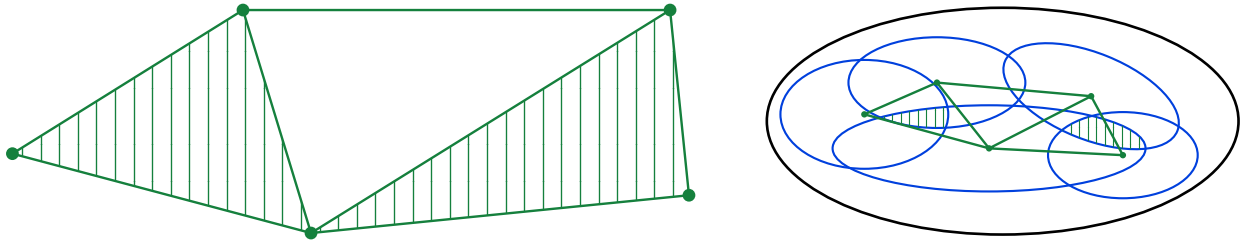

\centering
  \tfig5
  \vskip -.5pc
  \caption{The nerve of the cover, and the same simplicial complex drawn on $M$}\label{fig:5}
\end{figure}

First, Segal replaces the de Rham complex on the vertical axis with the
manifold~$M$.  And he replaces the \v Cech complex on the horizontal axis with
the \emph{nerve} of the cover.  The nerve is defined as follows.  Let $V$~be
the free vector space with basis~$A$.  For each $\alpha \in A$, the basis
vector~$p_\alpha \in V$ represents the open set~$U_\alpha $.  Start with the
disjoint union of the basis vectors.  For each pair $\alpha _0,\alpha _1\in A$,
if the intersection $U_{\alpha _0}\cap U_{\alpha _1}$ is nonempty, then adjoin
the 1-simplex (line segment) with vertices $p_{\alpha _0},p_{\alpha _1}$.  For
nonempty triple intersections $U_{\alpha _0}\cap U_{\alpha _1}\cap U_{\alpha
_2}$ adjoin the 2-simplex with vertices $p_{\alpha _0},p_{\alpha _1},p_{\alpha
_2}$, etc.  The nerve is illustrated in \autoref{fig:5}.
 
\begin{figure}[htbp]
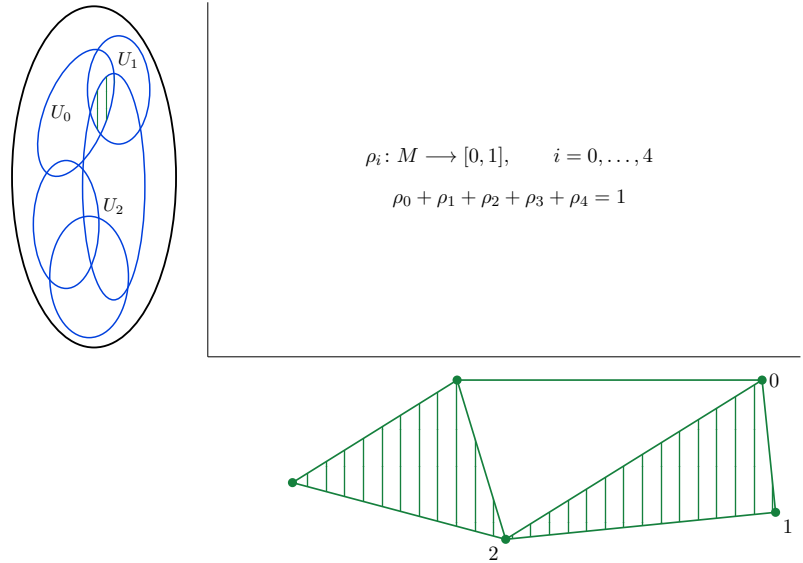

\centering
  \tfig6
  \vskip -.5pc
  \caption{A partition of unity as a function from~$M$ to the nerve}\label{fig:6}
  \end{figure}

A key ingredient in Weil's argument is a partition of unity subordinate to the
open cover: a collection of smooth functions~$\{\rho _\alpha \}_{\alpha \in A}$
on~$M$ with values in~$[0,1]$ such that at each~$m\in M$ only finitely
many~$\rho _\alpha $ are nonzero; $\operatorname{supp}\rho_\alpha\subset
U_\alpha$, so in particular if~$\rho _\alpha (m)\neq 0$ then $m\in U_\alpha $;
and $\sum\rho _\alpha (m)=1$.  In other words, the partition of unity
expresses~$m$ as a convex combination of finitely many open sets of the cover
in which it lies.  Segal's interpretation: the partition of unity is a function
from~$M$ to the nerve, as indicated in \autoref{fig:6}.

  \begin{figure}[htbp]
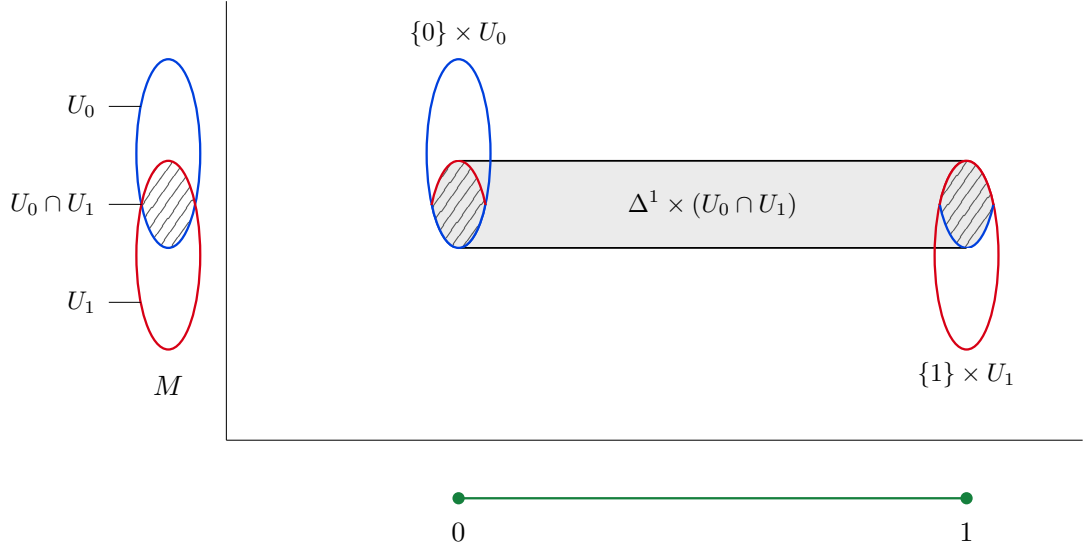

\centering
  \tfig7
  \vskip -.5pc
  \caption{Union of graphs of partitions of unity for two intersecting open sets}\label{fig:7}
  \end{figure}

\begin{figure}[htbp]
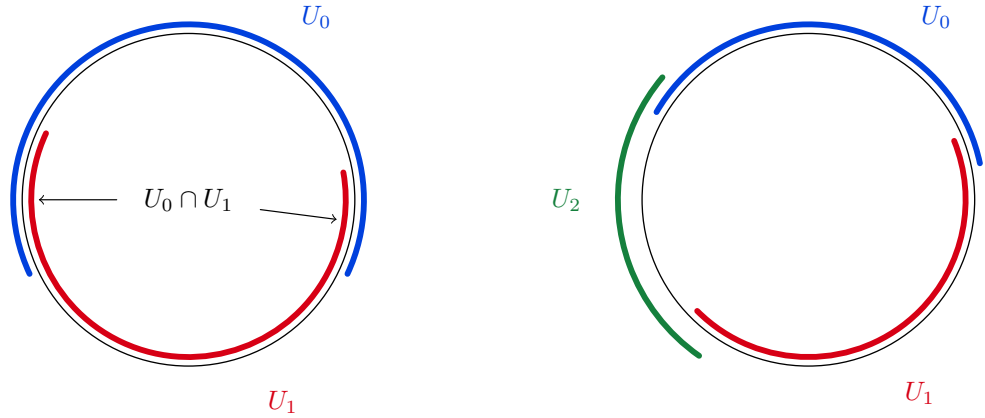

  \centering
  \tfig8
  \vskip +.5pc
  \caption{A cover of the circle by two arcs, and by three}\label{fig:8}
  \end{figure}

Next, Segal makes two maneuvers characteristic of a geometer: (1)~he replaces a
function by a \emph{space}, and (2)~rather than single out one particular
choice, he brings in the \emph{space} of all choices.  Applied to these
partitions of unity, he (1)~considers the graph of each partition of unity,
regarded as a map from~$M$ to the nerve of the cover, and (2)~he takes the
union of the graphs over all partitions of unity.  \autoref{fig:7} is the
result for the union of two intersecting open sets.  (The union of graphs~$W$
leads to a neat proof of the Mayer--Vietoris theorem.  The projection $W\to M$
is a homotopy equivalence, since its fibers are simplices, and for the subspace
$W'=\{0\}\times U_0 \;\cup \; \{1\}\times U_1$ excision gives $H^q(W,W')\cong
H^{q-1}(U_0\cap U_1)$.  Now look at the long exact sequence of the
pair~$(W,W')$.)  At this point I invite the reader to contemplate this
construction for the coverings of the circle indicated in \autoref{fig:8} as
well as a configuration of three open sets with a nonempty triple intersection.

Segal makes one more move to arrive at the classifying space of a category.
This was motivated by his wanting to also bring in Milnor's
construction~\cite{Mi1} of the classifying space of a topological group.  To
unify these constructions, Segal~\cite{S6} introduced the idea of a
\emph{topological category}~$C$.  Namely, $C$~is defined by \emph{topological
spaces}~$C^0,C^1$ of objects and morphisms, and the structure maps (such as the
maps $C^0 \leftleftarrows C^1$ that give the source and target of a morphism)
are required to be continuous.  A topological category~$C$ has a nerve, which
is a \emph{simplicial space}~$NC$:
  \begin{equation}\label{eq:2}
    \tsimpobj
  \end{equation}
Here $NC^0$~is the space of objects, $NC^1$~is the space of morphisms,
$NC^2$~is the space of composable ordered pairs of morphisms, etc.; the solid
maps delete an object, composing morphisms as needed, and the dotted arrows
insert an identity map.  A simplicial space has a \emph{realization} as a
topological space~\cite{Mi2}, all described elegantly in~\cite[\S1]{S6}.  The
realization of~$NC$ is called the \emph{classifying space}~$BC$ of the
topological category~$C$.

  \begin{figure}[htbp]
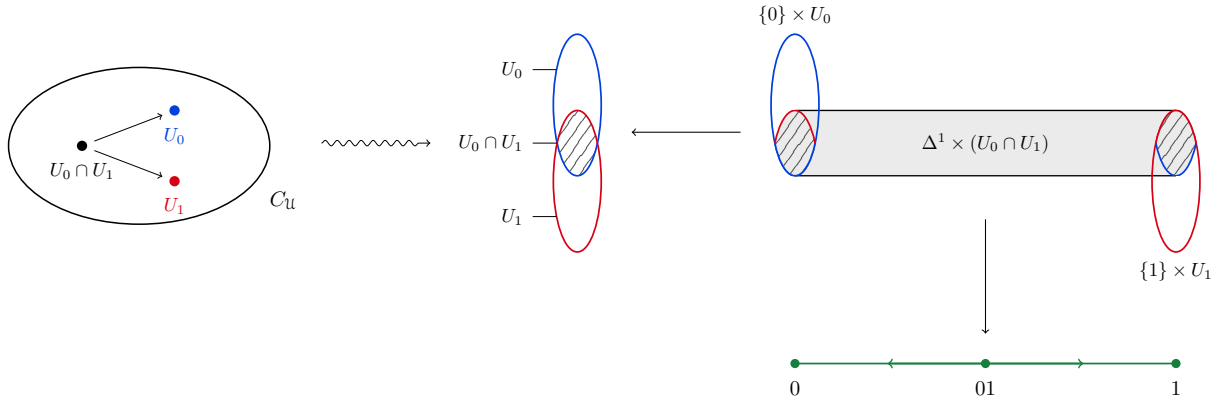

  \centering
  \tfig9
  \vskip -.5pc
  \caption{The topological category of a 2-set cover and its classifying space}\label{fig:9}
  \end{figure}

\autoref{fig:9} depicts the topological category~$C_{\sU}$ for a cover
$\sU=\{U_0,U_1\}$ consisting of two intersecting open sets.  The space of
objects~$(NC_{\sU})^0$ is the disjoint union of the spaces~$U_0$, $U_1$, and
$U_0\cap U_1$.  The morphisms in~$C_{\sU}$ map points of~$U_0\cap U_1$ to the
same point regarded as in~$U_0$ or~$U_1$; thus the space of non-identity
morphisms is the disjoint union of two copies of~$U_0\cap U_1$.  The
classifying space of this topological category is---up to homeomorphism---the
space we encountered in \autoref{fig:7}.  I remark in passing that
$NC_{\sU}$~is not merely a simplicial space; it is a simplicial
\emph{manifold}.

Now I return to Segal's motivating question ``\emph{How can I realize Weil's \v
Cech--de Rham double complex geometrically?}''  Weil's double complex is the de
Rham complex of the simplicial manifold whose space of $p$-simplices is
$\coprod U_{\alpha _0}\cap \cdots \cap U_{\alpha _p}$---a close cousin
of~$NC_{\sU}$ with the same realization up to homotopy---and its horizontal
differential is the alternating sum of the face maps.  The double complex
computes the cohomology of the classifying space~$BC_{\sU}$.  The two
projections of~$BC_{\sU}$---to~$M$, with simplices as fibers, and to the
nerve---lead to two spectral sequences of the double complex.  Spectral
sequences are the subject of~\cite[\S5]{S6}.

Segal's short, elegant paper~\cite{S6} was an aside on the way to his
thesis~\cite{S12}.  It had an immediate application to the filtration in
$K$-theory, as explained in~\cite[\S5]{S6}; previously the filtration depended
on a cellular structure.  There was some skepticism about the significance of
the general picture, but fortunately Grothendieck encouraged Segal to publish
his paper.  Over time it became a landmark.  Its influence has been vast and
long-lasting; the construction is firmly embedded into modern algebraic
topology.  The account I gave emphasizes the small-scale elegance and thought
process.  But---in a geometric sense---the significance of its output is large
scale: the \emph{homotopy type} of the classifying space.  Constructing
homotopy types by geometric means is a theme that runs through Segal's work.
 
If a topological category~$C$ possesses extra structure, then naturally so too
does its classifying space~$BC$.  In a subsequent paper~\cite{S3}, Segal
identifies a structure on~$C$ that leads to an \emph{infinite loop structure}
on~$BC$.  This arose from two observations.  First, Segal characterized a
symmetric composition law on a category, such as direct sum of vector spaces
or disjoint union of finite sets, in terms of its universal property.  So, if $X_1,X_2$~are objects
of~$C$, then a composition~$Y$ is an object equipped with maps $X_1\to Y\gets
X_2$ with the universal property that in any diagram of the following shape a
unique dotted morphism exists: 
  \begin{equation}\label{eq:3}
     \begin{gathered}
     \xymatrix@R-1.5pc@C+2pc{&Y\ar@{-->}[dd]\\X_1\ar[ur]\ar[dr]&&X_2\ar[ul]\ar[dl]\\
     &Y'} \end{gathered}
  \end{equation}
Therefore, if $C_2$~is the category whose objects are composition diagrams
$X_1\to Y\gets X_2$, then in the diagram of forgetful maps
  \begin{equation}\label{eq:4}
     \begin{gathered} \xymatrix@R-.7pc{&C_2\ar[dl]\ar[dr] \\ C\times C&&C}
     \end{gathered} 
  \end{equation}
the left arrow, mapping $X_1\to Y\gets X_2$ to~$(X_1,X_2)$, is an equivalence
of categories.  The correspondence diagram~\eqref{eq:4} expresses the
composition law.  Furthermore, transposition of $X_1,X_2$ acts on~$C_2$; the
left arrow is equivariant and the right-hand arrow is strictly invariant.
Commutativity is built in, with no coherence conditions to check.  
 
A related second observation concerns the homotopy category of topological
spaces, in which the objects are topological spaces and morphisms are homotopy
classes of continuous maps.  Namely, a symmetric map $X\times X\to Y$ in that
category is a diagram
  \begin{equation}\label{eq:5}
     \begin{gathered} \xymatrix@R-.7pc{&X'\ar[dl]\ar[dr] \\ X\times X&&Y}
     \end{gathered} 
  \end{equation}
in which $X'$~is a topological space equipped with a $\Sym_2$-action, the left
arrow is a $\Sym_2$-equivariant homotopy equivalence, and the right arrow is
$\Sym_2$-invariant.

These observations led Segal to introduce $\Gamma $-categories and $\Gamma
$-spaces.  Iterating the first, one obtains for each positive integer~$k$ a
category~$C_k$ of $k$-fold composition diagrams, each equivalent to~$C^{\times
k}$, and functorial in the finite pointed set $\{0,1,\dots ,k\}$.  This is a
\emph{$\Gamma $-category}; applying the classifying space construction gives a
\emph{$\Gamma $-space}: a functor~$\Phi$ from finite pointed sets to spaces with
$\Phi(\{0\})$~contractible and the \emph{Segal maps} $\Phi(\{0,\dots ,k\})\to
\Phi(\{0,1\})^{\times k}$ homotopy equivalences.  That last requirement---the
\emph{Segal condition}---is the analog of the left arrow in~\eqref{eq:5} being
a homotopy equivalence.  Segal used $\Gamma $-spaces to construct a ``delooping
machine''.  These ideas not only found immediate applications at the time they
were introduced, but they are also fundamental in current work, for example in
the theory of higher categories.  Here is a sample of important work that
builds on~\cite{S3}, a somewhat arbitrary selection from many possibilities:
\cite{BF, MT, T, R1, B, L, R2, Si, GGN, H, CS, DK, BSP}.

   \section{Quantum field theory}\label{sec:4}
% lastsubsec@000

\begin{figure}[ht]
  \centering
  \includegraphics[scale=.35]{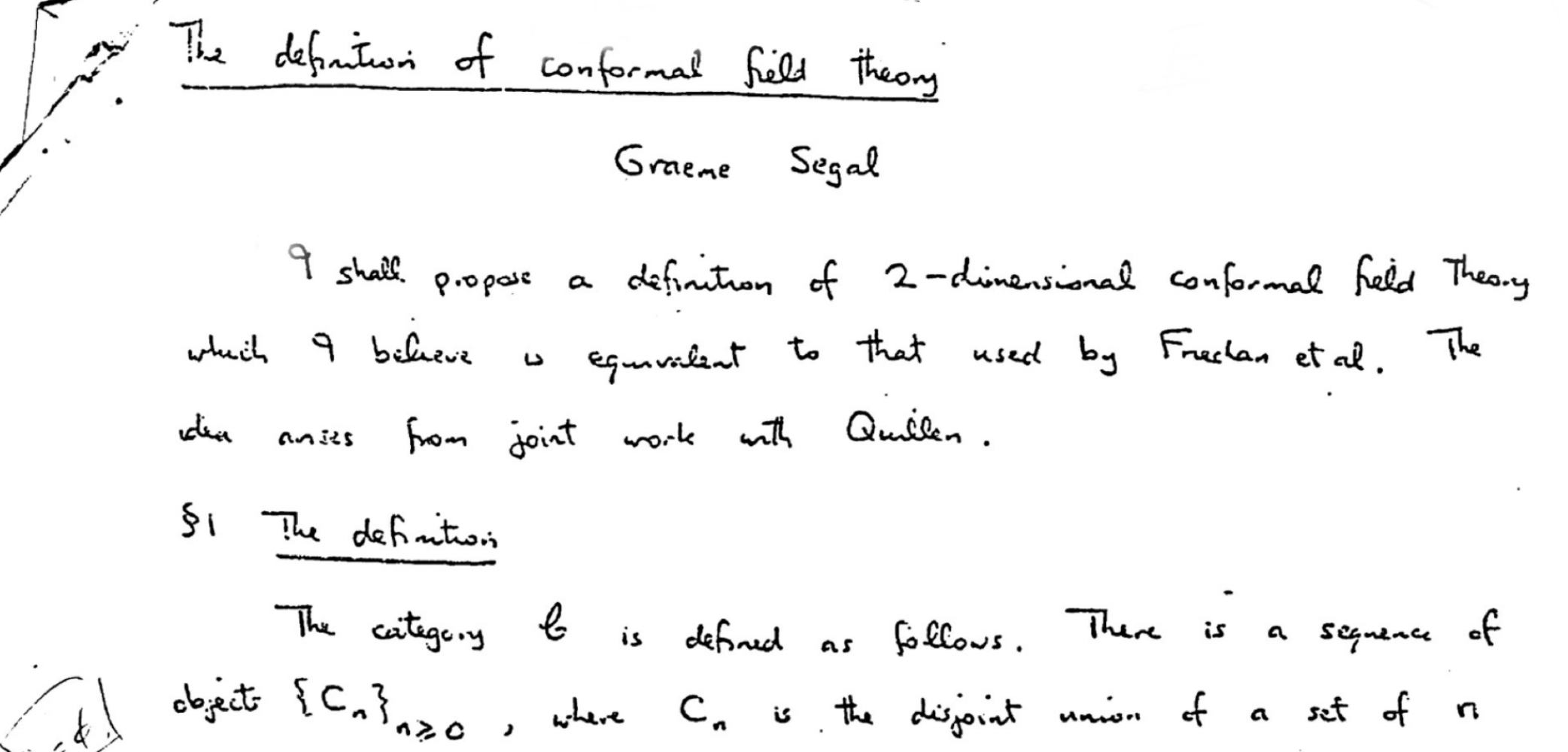}
  \vskip -.5pc
  \caption{Segal's 1987 preprint}\label{fig:10}
  \end{figure}

In 1987 Segal introduced a new axiom system for 2-dimensional conformal field
theory.  At the time, such theories were under intense investigation by
theoretical physicists, due to their role in string theory.  The work
interacted deeply with Riemann surface theory on the one hand and with the
representation theory of infinite-dimensional Lie groups and Lie algebras on
the other.  Segal had been working on loop groups and their
representations~\cite{S5,PS}, and he also had a longstanding interest in
quantum field theory~(QFT), dating back to his graduate student days.  So it
was natural that he engage with this new phase of the math--physics
interaction.  His ideas were written in a preprint (\autoref{fig:10}) that
circulated for a while as a samizdat manuscript, though a very brief account
appeared in~\cite{S19}; the full published version~\cite{S10} only appeared
after the turn of the millennium.

Here, in Segal's words, is his journey from quantum field theory to geometry
and topology as a graduate student in the 1960s, and then back to quantum field
theory in the 1980s:

  \smallskip
  \begin{quotation}\footnotesize
 I arrived in Cambridge as a graduate student in 1962, intending to study
mathematical physics.  But, in the days before the Standard Model was invented,
I was bewildered and discouraged by the lectures on quantum field theory I went
to, with their focus on the ``analytic $S$-matrix'' then in vogue.  Instead, I
was attracted to topology and geometry, and felt very proud when Hodge agreed
to be my supervisor [\dots]

I'd begun trying to learn QFT more systematically around 1980.  Two-dimensional
spacetime seemed the easiest starting-point, so I studied the Glimm--Jaffe
path-integral construction of 2-dimensional massive bosonic theory with an
arbitrary polynomial interaction.  They constructed a theory for each
polynomial, but I was struck that there seemed no way of knowing whether these
theories were actually different, and still less of recovering the interaction
polynomial from the theory.  I also tried to understand the claimed
equivalence, in two space-time dimensions, between bosonic theory with
sine-Gordon interaction and fermionic theory with a quartic Thirring
interaction.  Renormalization still bewildered me, but during a term in
Berkeley in 1982 I learned about Friedan's work relating the renormalization
group flow of a 2-dimensional $\sigma$-model to the Ricci curvature of its
target.  I slowly assimilated the idea that field theories are defined by
Lagrangians which serve as the weighting factors of path-integrals, and,
crucially, that the renormalizability of the theory is equivalent to the
well-posedness of the variational problem defined by the Lagrangian.  This
suggested that a geometric approach to the subject was possible, but there
seemed an urgent need for a better mathematical definition of a quantum field
theory---one which would make it possible to speak of spaces of theories of
various kinds.

The final stimulus before I had anything to say of my own was learning from
Witten---and the active group of younger people around him---about the role of
2-dimensional conformal field theory in string theory.  The central idea of
string theory (at that time!)\ was that the space of Riemannian manifolds is
not a satisfactory configuration space for quantum gravity, and should be
replaced by the space of 2-dimensional QFTs.  A Riemannian manifold~$M$ gives a
2-dimensional field theory---its $\sigma$-model---in which the fields on a
space-time~$\Sigma$ are the maps from~$\Sigma$ to~$M$.  Friedan's thesis showed
that Ricci flow on Riemannian manifolds corresponds to the renormalization flow
on 2-dimensional QFTs, and Zamolodchikov's `$c$-theorem' asserted that the
renormalization flow of 2-dimensional QFTs is along the gradient of a natural
real-valued function~$c$ on the space of theories.  This looked close to being
a precise mathematical theorem, and became my lodestar.

Talking to Witten, I reformulated---at first as a joke---his account of QFT in
his plenary lecture~\cite{Wi1} at the 1986 ICM as defining a field theory in
$d$-dimensional space-time as a functor from the $d$-dimensional cobordism
category to the category of topological vector spaces.  In other words, it is a
representation of a category, and so it is a structure not so different from a
representation of a Lie group or semigroup.  In the special case of
2-dimensional conformal field theories it proved surprisingly easy to develop
this picture in detail, and the account I wrote at the end of the 1980s was
circulated widely (not by me), and eventually published.
  \end{quotation}
  \smallskip

\noindent
 Today we use Segal's axiom system for general quantum field theories, though
as he writes ``it took some time before I felt confident that the cobordism
definition could be applied beyond its original application to conformally
invariant theories.''
 
My account here emphasizes the formal mathematical structure and mathematical
applications.  It leaves out Segal's physical and philosophical intuitions that
inspired the axioms, and also does not delve into the impact back in physics.
I begin by recalling older axiom systems for quantum theory, placing Segal's
axiom system in a broader context.  Then I introduce his axioms, discuss some
of his work in this area, and describe a few of the mathematical ramifications
of Segal's conceptualization of quantum field theory.  As I write there is
still much to do.  For more detailed explanations of the axiom systems
described here, including Segal's, see~\cite{F}.
 
In the early days of quantum mechanics, around 1930, Dirac~\cite{D} and von
Neumann~\cite{vN1,vN2} gave the first axiom system for quantum theory; it was
nicely codified later by Mackey~\cite{Mac}.  From a complex separable Hilbert
space~$\sH$ they derive: a convex space of \emph{states}, a $*$-algebra of
\emph{observables}, and a pairing between them to probability measures
on~$\RR$.  The extreme states, called pure states, form the projective
space~$\PP\sH$.  The $*$-algebra of bounded operators on~$\sH$ is the algebra
of observables; the self-adjoint operators pair with states.  There is also
dynamics, a strongly continuous unitary representation 
  \begin{equation}\label{eq:6}
     U\:\RR\longrightarrow  U(\sH) 
  \end{equation}
of the additive group of time translations.  This axiom system puts states and
observables on an equal footing.  An alternative axiom system by Irving
Segal~\cite{Se} in the 1940s favors observables as the primary concept, leaning
on the then emerging field of operator algebras.  One can loosely think of
these two approaches as the Schr\"odinger/Heisenberg dichotomy between
emphasizing states and emphasizing observables.
 
The relativistic setting for quantum field theory is Minkowski spacetime, an
affine space~$M$ over a real vector space equipped with a Lorentz signature
nondegenerate bilinear form and a time-orientation.  (The one-dimensional case
reduces to quantum mechanics.)  Here too there are two commonly used axiom
systems, dating from the late 1950s and early 1960s.  The Haag--Kastler
axioms~\cite{HK} fall on the Heisenberg side of the ledger in that they
emphasize observables over states.  But now there is a \emph{local} theory of
observables, encoded in a net of operator algebras that associates an algebra
of operators to regions of~$M$.  The Wightman axioms~\cite{StW,Ka} use fields
as observables and include a vacuum state.  A unitary representation of the
Poincar\'e group generalizes the dynamics in the Dirac--von Neumann axioms.
There is a nicely developed theory here, especially for the Wick-rotated
version in Euclidean space and theorems~\cite{OS2} that mediate between the
Euclidean and Minkowskian variants.
 
\begin{figure}[ht]
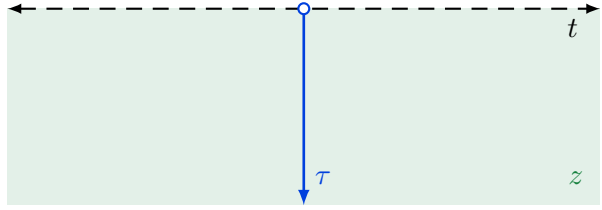

\centering
  \tfig{10}
  \vskip -.5pc
  \caption{Wick rotation in quantum mechanics}\label{fig:11}
\end{figure}

(Graeme) Segal's axioms for QFT are in the Wick-rotated setting, so I pause to
introduce Wick rotation in quantum mechanics.  By Stone's theorem there is a
self-adjoint operator~$H$ on~$\sH$, generally unbounded, that generates the
unitary representation~\eqref{eq:6} in the sense that $U_t = e^{-itH/\hbar}$.
Introduce the crucial assumption that the \emph{Hamiltonian}~$H$ is bounded
below, and by shifting by a constant we may assume that $H$~is nonnegative.
Then the unitary time evolution $t\mapsto e^{-itH/\hbar}$ is the boundary value
of a holomorphic one-parameter semigroup $z\mapsto e^{-izH/\hbar}$ of
contractive operators; here $z\in \CC$ satisfies $\Im(z)<0$.  Wick rotation is
the restriction to $z=-i\tau $, $\tau \in \Rp$, which is a real one-parameter
semigroup $\tau \mapsto e^{-\tau H/\hbar}$ of contractive operators; see
\autoref{fig:11}.  For example, in the case of a quantum mechanical particle of
mass~$m$ moving on a closed Riemannian manifold~$N$, the Hilbert space~$\sH$ is
the space of $L^2$~complex-valued functions on~$N$; the Hamiltonian is
$\hbar^2\Delta /2m$, where $\Delta $~is the Hodge laplacian on functions; and
the real semigroup $\tau \mapsto e^{-\tau \hbar\Delta /2m}$ solves the heat
equation on~$N$.  These heat operators are analytically more tractable than the
unitary operators defined from the laplacian for real~$t$; for example, they
are smoothing operators.  The relative analytic simplicity makes the
Wick-rotated setting more attractive for mathematical analysis.
 
\begin{figure}[htbp]
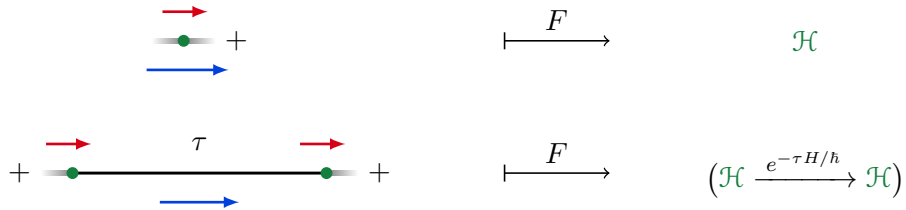

\centering
  \tfig{11}
  \vskip -.5pc
  \caption{A representation of oriented Riemannian closed bounded
  intervals}\label{fig:12}
\end{figure}

One route to motivate Segal's axioms is to write the real semigroup $\tau
\mapsto e^{-\tau H/\hbar}$ in a geometric form.  Namely, Riemannian manifolds
diffeomorphic to the closed bounded interval~$[0,1]$ form a real semigroup
under concatenation.  The only invariant is the total length~$\tau $, and so
this semigroup (of equivalence classes) is isomorphic to~$\Rp$.  The
illustration in \autoref{fig:12} represents these intervals as \emph{bordisms}
from a positively oriented point to a positively oriented point.  The red
arrows indicate which boundary components are incoming in the bordism and which
are outgoing.  The blue arrows indicate the 1-dimensional orientation.  The
point itself is mapped to the Hilbert space~$\sH$, and the bordism to the
contractive map $e^{-\tau H/\hbar}$.  The gray bands around the points are
embeddings into a germ of a Riemannian 1-manifold.  These germs allow us to
glue the 1-dimensional bordisms along open intervals, which smoothly patches
manifolds and metrics.  The red arrow over the point is then a coorientation
inside the germ; it is the Wick-rotated remnant of the time orientation in
Minkowski spacetime.
 
\begin{figure}[ht]
\centering
  \tfig{13}
  \vskip -.5pc
  \caption{Wick-rotated QFT as a homomorphism $\tbord(\sF)\to \tVect$}\label{fig:13}
\end{figure}

Segal's axioms assert that a Wick-rotated QFT is a homomorphism of symmetric
monoidal categories,\footnote{Just as the semigroup~$\Rp$ lacks an identity
element, so too these categories lack identity morphisms.} from a bordism
category to a linear category of topological vector spaces.  As with all great
mathematical conceptualizations, a few words pack lots of punch.  Let us
unpack.  There are two parameters that define the domain bordism category:
$n,\sF$.  The positive integer~$n$ is the dimension of ``spacetime'', though
that nomenclature refers to the QFT before Wick rotation.  I come to~$\sF$
presently.  The bordism category $\tbord$ has as objects closed
$(n-1)$-manifolds~$Y$.  They are embedded in a germ of an $n$-manifold and are
cooriented in that germ.  Morphisms are (isomorphism classes of) bordisms:
compact $n$-manifolds with boundary, a partition of the boundary into two
unions of components (``incoming'' and ``outgoing''), and identifications of
these unions with the domain and codomain of the morphism.  The symmetric
monoidal structure is disjoint union.  The codomain linear category in Segal's
axioms is a category of topological vector spaces and nuclear maps between
them.  Segal's axiom system is the following.

  \begin{definition}[]\label{thm:1}
 A \emph{field theory} is a continuous homomorphism 
  \begin{equation}\label{eq:7}
     F\:\tbord(\sF)\longrightarrow \tVect 
  \end{equation}
of symmetric monoidal topological categories.
  \end{definition}

\noindent
 This formulation is not quite the one in the first paper~\cite{S10}.  It is
one that developed over time and appears in the joint paper~\cite{KS} with
Kontsevich.  
 
In the case of quantum mechanics (\autoref{fig:12}) the 1-dimensional intervals
carry an orientation and a Riemannian metric.  One-dimensional orientations and
Riemannian metrics are examples of \emph{local} structures on 1-manifolds that
pull back under local diffeomorphisms.  For a general field theory the
\emph{sheaf}~$\sF$ carries these structures.  The pair~$(n,\sF)$ determines the
bordism category and specifies the \emph{type} of the field theory.  Let
$\Man_n$ be the category whose objects are smooth $n$-manifolds and whose
morphisms are local diffeomorphisms.  It carries a Grothendieck topology
defined as usual by open covers.  The \emph{fields} of a Wick-rotated field
theory form a sheaf~$\sF$ on~$\Man_n$ with values in the category of sets or
simplicial sets.  In the bordism category $\tbord(\sF)$ all objects and
morphisms are equipped with sections of the sheaf~$\sF$ (over the germs).  The
theory is \emph{topological} if it factors through a bordism category in which
this sheaf is locally constant, in which case it is equivalent to a tangential
structure in the sense of classical bordism theory~\cite{La}.

  \begin{remark}[]\label{thm:2}
 Topological bordism categories were introduced by Milnor~\cite[\S1]{Mi3} in
his proof of the $h$-cobordism theorem.
  \end{remark}

  \begin{remark}[]\label{thm:3}
 In terms of the Schr\"odinger/Heisenberg dichotomy I evoked earlier, Segal's
axioms fall into the Schr\"odinger picture.  A set of Heisenberg-like axioms
for Wick-rotated field theory---a Wick-rotated version of Haag--Kastler in the
spirit of Segal's axioms---is developed by Costello--Gwilliam
in~\cite{CG1,CG2}.
  \end{remark}

It would be difficult to overestimate the profound effect of Segal's axioms.
They gave many mathematicians a foothold in quantum field theory, especially on
the more geometric side that had been growing in the hands of Witten and other
physicists.  There were immediate applications of 2-dimensional conformal field
theory already in~\cite{S10}, for example to the modular behavior of positive
energy representations of loop groups.  One outcome of discussions between
Segal, Atiyah, and Witten was Witten's formulation~\cite{Wi2} of the Jones
invariants of knots and links in terms of quantum field theory.  This in turn
led to Atiyah's axiomatization~\cite{A} of topological field theory and the
emphasis on the relationship to Thom's bordism theory.  Many variations were
developed and applications made.  Over the almost forty years since their
introduction, Segal's axioms catalyzed great progress on the
geometric/topological/algebraic side of QFT.  Of course, there were analytic
developments as well during this period, and we are beginning to see more and
more crosstalk between these different strands in the mathematical conception
of QFT.
 
Segal has contributed to this endeavor in several works subsequent
to~\cite{S10}.  His six Stanford lectures~\cite{S13--S17} in the mid-1990s were
never completely written up, but still had wide influence.  For example, in
those lectures Segal introduces cochain-valued field theories and even
spectrum-valued field theories, which stimulated a homotopy-theoretic take on
this mathematical structure.  He also indicated some relationships---developed
later---between Atiyah--Singer index theory and field theory.  In addition, he
gave first ideas of extended locality, which he called a \emph{3-tier theory}.
Segal returned to several of these topics in subsequent work.

One example is the work~\cite{CJS} on Floer homotopy theory, as already
mentioned in~\S\ref{sec:1}, which picks up on the homotopy-theoretic side of
field theory.  It is closely related to the quantization of fermionic fields
and to semi-infinite cohomology, both topics that had occupied Segal in various
contexts.  Moreover, it fits the theme in Segal's work of attaching homotopy
types to various mathematical and physical structures.
 
Segal took up extended locality in subsequent papers.  In joint work~\cite{MoS}
with Greg Moore, Segal considers 2-dimensional topological field theories.  It
was well-known, since Dijkgraaf's thesis~\cite{Dijk}, that a 2-dimensional
topological theory of oriented manifolds is given by a commutative Frobenius
algebra.  Moore--Segal extend this to an ``open-closed theory'', closely
related to a 3-tier extension.  In a 3-tier theory a point is assigned a linear
category; here that role is played by the category of boundary conditions
(``D-branes'').  In fact, the focus is not on the 3-tier extension, but rather
on the interpretation as the baby case of a string theory background.
 
The 3-tier extension is in the foreground in~\cite{S18}.  In an $n$-dimensional
Wick-rotated QFT a closed $(n-1)$-manifold~$Y$ maps to a topological vector
space.  Here Segal considers a decomposition $Y=Y_1\cup \mstrut _ZY_2$ as the
union of two $(n-1)$-manifolds with common boundary a codimension-one
submanifold~$Z\subset Y$.  In topological field theory such decompositions are
common: one may assign an algebra~$A$ to~$Z$, left and right $A$-modules
to~$Y_1,Y_2$, and then the vector space attached to~$Y$ is the tensor product
over~$A$ of these modules.  Here Segal considers nontopological theories, such
as the free spinor field.  Then one must thicken~$Z$ to an overlapping region
in the intersection $Y_1\cap Y_2$, and the attached algebra is a von Neumann
algebra.  The tensor product is now the more intricate one introduced by
Connes.  
 
In the mathematical developments surrounding relativistic quantum field theory,
the important Osterwalder--Schrader reconstruction theorem~\cite{OS1,Gl,OS2,GJ}
produces the Minkowski spacetime theory from the Euclidean theory.  The joint
paper~\cite{KS} with Kontsevich tackles the reconstruction question on curved
manifolds.  They introduce what we might think of as a geometric version of the
holomorphic semigroup depicted in \autoref{fig:11}, now in arbitrary dimension.
This is the notion of an \emph{allowable complex metric} on a smooth manifold.
The space of allowable complex metrics on a vector space is a contractible
domain that contains the Riemannian metrics, and the nondegenerate points of
its Shilov boundary are the time-oriented Lorentz metrics.  That is an accurate
picture on a single vector space, but over a smooth manifold whereas allowable
complex metrics always exist (they are sections of a fiber bundle with
contractible fibers), there are global topological obstructions to the
existence of Lorentz metrics.  Kontsevich--Segal prove that a unitary theory on
allowable complex metrics induces a functor on compact, real-analytic, globally
hyperbolic Lorentzian bordisms, which assigns Hilbert spaces to time-symmetric
hypersurfaces and unitary operators to bordisms between them.  Their paper
contains many important insights beyond this reconstruction theorem and will be
a springboard for future work.

Finally, I mention that Segal has recently made progress in enhancing his axiom
system to encode good small distance behavior of the theory.  I refer to his
contribution to these proceedings.
 
To conclude, I return to the scales of \autoref{fig:3}.  Segal's axiom system
is large-scale mathematics at its finest: a conceptualization of structure
whose impact extends across mathematics and theoretical physics, and which
continues to open new directions.

%\appendix

 \bigskip\bigskip
\providecommand{\bysame}{\leavevmode\hbox to3em{\hrulefill}\thinspace}
\providecommand{\MR}{\relax\ifhmode\unskip\space\fi MR }
% \MRhref is called by the amsart/book/proc definition of \MR.
\providecommand{\MRhref}[2]{%
  \href{https://mathscinet.ams.org/mathscinet-getitem?mr=#1}{#2}
}
\providecommand{\href}[2]{#2}

  \end{document}